\documentclass[12pt,a4paper,xdvi,oneside]{article}
\usepackage{epsfig}
\usepackage{graphicx}

\usepackage{float}
\usepackage{paralist}
\usepackage{lscape}
\usepackage{comment}

\usepackage{ae,aecompl}   
\usepackage{amsfonts,amssymb,amsthm}
\usepackage{latexsym}
\usepackage{amsmath}
\usepackage{times}

\begin{document}

\newtheorem{definition}{Definition}
\newtheorem{lemma}{Lemma}
\newtheorem{corollary}{Corollary}
\newtheorem{theorem}{Theorem}
\newtheorem{hypothesis}{Conjecture}
\newtheorem{proposition}{Proposition}
\newtheorem{remark}{Remark}
\newtheorem{example}{Example}

\newcommand{\settocdepth}[1]{\addtocontents{toc}{\protect\setcounter{tocdepth}{#1}}}
\newcommand{\Prob}{\mbox{\bf P}}
\newcommand{\calA}{\mbox{${\cal A}$}}
\newcommand{\calC}{\mbox{$\cal C$}}
\newcommand{\equivA}{\mbox{$\equiv_{\cal A}\,$}}
\newcommand{\Transp}{\,\!^{\mbox{\footnotesize\bf T}}}

\newcommand{\Rspace}{\mbox{${\mathbb R}^N$}}
\newcommand{\Mspace}{\mbox{${\mathbb M}(N\times N,{\mathbb R})$}}
\newcommand{\maV}{\mbox{${\mathbb V}$}}
\newcommand{\maM}{\mbox{${\mathbb M}$}}
\newcommand{\SSy}{\mbox{$StatSym$}}

\newcommand{\SyG}{\mbox{${\cal S}$}}
\newcommand{\AlG}{\mbox{${\cal A}$}}
\newcommand{\CyG}{\mbox{${\cal Z}$}}
\newcommand{\GG}{\mbox{${\cal G}$}}
\newcommand{\sGG}{\mbox{$\cal \!\;_G$}}

\restylefloat{table}
\restylefloat{figure}
\setcounter{secnumdepth}{3}
\setcounter{tocdepth}{2}

\setcounter{page}{0}

\title{\vspace{-5mm} State Symmetries in Matrices and Vectors\\
on Finite State Spaces}

\author{Arne Ring\thanks{arne.ring@epost.de} \\ 
D-88397 Biberach/Ri\ss
}

\maketitle
\thispagestyle{empty}

\section*{Abstract}

State symmetries are defined as permutations which act on vector spaces of \mbox{column} vectors 
and square matrices, resulting in isotropy groups for specific \mbox{vector} spaces.

A large number of properties for such objects is shown, to provide a rigorous basis for 
future applications.
The main statement characterises the state symmetry of vector sequences $(v^{(i)})$ which are generated by 
powers of a generator matrix~$M$: \mbox{$v^{(i)}= M^i\,v^{(0)}$}.

A section of examples illustrates some applications of the theory.

\subsection*{Keywords}
state symmetry, linear algebra, vector space, symmetric group, general linear group, isotropy group, vector sequences.

\subsection*{Mathematics Subject Classification (MSC 2000)}
05E20 Group actions on designs, geometries and codes, \\
15A03 Vector spaces, linear dependence, rank,\\
20G20 Linear algebraic groups over the reals, the complexes, the quaternions.

\newpage

\tableofcontents

\subsection*{Publishing Notification}
The results of this paper had been presented at the 4ecm 
(4th european conference of mathematics in Stockholm, June 2004) \cite{Ring:04a}.

\newpage


\section{Introduction}

This paper origins in the analysis of Markov chains which act on finite state spaces with 
"some symmetries" \cite{Ring:96}.
If the transition graphs of the Markov chain has a symmetric structure, it may be lumped 
to a chain with a smaller state space -- which may reduce the calculation expense. 
In addition, new properties of Markov chains and their limit distributions can be derived \cite{Ring:02}.

The investigations led to the finding that the underlying theory of transition matrices and start vectors 
can be reduced to problems of linear algebra and group \mbox{theory}. 
This theoretical part was separated from the application to Markow chains and shall now be presented 
in definitions and properties of special matrices and vectors. 
Theorems and conclusions will then serve as a basis for applications, 
not only restricted to Markow chains but also applicable e.g. to graph theory.
\medskip

The usual matrix symmetries is defined on the basis of $M\Transp=M$, where $M\Transp$ is the transposed matrix. 
In some cases they help to reduce the calculation expense; hence they may be used e.g. in numerical algorithms.

However, our term ''state symmetry'' is different to the usual matrix symmetry, although it will also be defined 
for square matrices (and for column vectors). It is based on transformations of permutation matrices
(which is called ''reordering'') and investigates the invariant objects of this transformation.

With other words, state symmetries are defined on an underlying finite ''state space'' $C$. 
Let $|C|=N$ be the number of states, then the $N\times N$ matrix $M$ 
is defining a weigthed graph~$G$ of edges between the vertices in $C$.
We will see that state symmetries of $M$ are based on special reorderings of entries the matrix $M$, 
which, roughly \mbox{speaking}, keep the structure of the graph $G$.

\medskip

State symmetries of matrices and vectors lead to two distinct structures: 
Once, the permutations of the state space, forming subgroups of the symmetric group~$S_N$; 
second, state symmetric matrices and vectors forming subspaces of vector spaces~$\Rspace$ and
$\Mspace$ and subgroups of the general linear group~$GL_N({\mathbb R})$ 
(we are using the field $\mathbb R$ of real numbers, although the statements 
are also valid for any other field).

The following paper presents the properties of both structures: 
After the formal definitions of reordering and state symmetries as application of permutations, 
we show properties of permutations which belong to given matrices and vectors (Section~2).
In Section~3, properties of state symmetric vectors and matrices are shown which arise for a given permutation of $C$.
The main statement characterises the state symmetry of vector sequences $(v^{(i)})$ which are 
generated by powers of a~generator matrix $M$: \mbox{$v^{(i)}= M^i\,v^{(0)}$}.

\medskip

A section of applications and a section with examples are 
added to illustrate the results of the theory.

\pagebreak  
\section{State symmetries and their structure}
\subsection{Permutations of the state space}

Let $C=\{1,2,\ldots,N\}$ be a finite set of natural numbers, which is called the state space. 
Permutations $\rho$ are elements of the symmetric group $\SyG_N$, and they may act on the state space $C$ so that 
its elements are changed with each other ($\rho: c \mapsto \rho(c)$).
It is now natural to ask about e.g. invariants of this mapping (stabilizers) and on orbits of a permutation 
within a state space.

In this section we want to apply the action of the symmetric group $\SyG_N$ on vector spaces 
of column vectors $\Rspace$ and of square matrices $\Mspace$ of real numbers. 

\subsection{Reordering} \label{reordering}

We start with the definition of reordering:

\begin{definition}[Reordered matrices and vectors] \label{defreorder} ~
Let $C=\{1,\ldots,N\}$ be a state space with $N<\infty$ elements and let $\rho$ be a permutation of $C$.
Let $v=\left(v_{i}\right)$ be a column vector with dimension $N$ and 
$M=\left(m_{i\,j}\right)$ a $N\times N$ square matrix of real numbers.
\begin{enumerate}
\item $v^\rho=\left(v_i^\rho\right):=\left(v_{\rho(i)}\right)$ is called the \underline{reordered} vector.
\item \label{defreordmat}
	$M^\rho=\left(m_{i\,j}^\rho\right):=\left(m_{\rho(i)\,\rho(j)}\right)$ is called the \underline{reordered} matrix. 
\end{enumerate}
\end{definition}

Reordering defines an action of elements of the symmetric group $\SyG_N$ on the sets $\Rspace$ and 
$\Mspace$: For the identity permutation $id$ it holds \mbox{$M^{id}=M$} and $v^{id}=v$,
and for all $\rho_1,\rho_2$, $M$ and $v$ it holds 
$M^{\rho_1\circ\rho_2}= (M^{\rho_2})^{\rho_1}$ and \mbox{$v^{\rho_1\circ\rho_2}= (v^{\rho_2})^{\rho_1}$}.\\
Hence, $\Rspace$ and $\Mspace$ are $\SyG_N$-sets.
\medskip

For an alternative way of writing the reordering, linear transformations may be used.
If $\rho$ is a permutation, the corresponding permutation matrix $T=(t_{ij})$ is given by
\begin{equation}
t_{ij}=
  \begin{cases}
    1 & \mbox{if~~} j=\rho(i)  \\ 
    0 & \mbox{otherwise .}
  \end{cases}
\end{equation}
Then, the reordered matrix is given by $M^\rho=T\,M\,T^{-1}$, and the reordered vector by $v^\rho=T\,v$.
As a result, $M$ and $M^\rho$ are similar matrices.
Note that permutation matrices are orthogonal ($T^{-1}=T\Transp$).
\smallskip

Reordering of objects leads to similar objects, which have the same elements but at different places.
This leads to special properties, like invariantness with regard to typical matrix calculations. 

\begin{lemma}\label{lemreorder}
Let $C$ be a finite state space, $\rho$ be a permutation of $C$ (with the \mbox{permutation} matrix $T$), 
$M$, $M_1$ and $M_2$ be $N\times N$ matrices and $v$ and $w$ be $N$-dimensional vectors. \\
Then the following properties hold:
\begin{enumerate}
\item The transposed reordered matrix is equal to the reordered transposed matrix 
$(M^\rho)\Transp =(M\Transp)^\rho$.
\item If $M$ is invertible, then  $(M^\rho)^{-1} =(M^{-1})^\rho$ .
\item Addition: $v^\rho+w^\rho = (v+w)^\rho$, ~ $M_1^\rho+M_2^\rho = (M_1+M_2)^\rho$ .
\item Multiplication: $M_1^\rho\,M_2^\rho = (M_1\,M_2)^\rho$, ~$M^\rho\, v^\rho=(M\, v)^\rho$ .
\item $\lambda\,M^\rho = (\lambda\,M)^\rho$, ~$\lambda\, v^\rho=(\lambda\, v)^\rho$.
\end{enumerate}
\end{lemma}
\begin{proof}
\begin{enumerate}
\item $(M^\rho)\Transp = (T\,M\,T^{-1})\Transp = (T^{-1})\Transp\,M\Transp\,T\Transp=T\,M\Transp\,T^{-1}=(M\Transp)^\rho$,
because $T$ is orthogonal.
\item $(M^\rho)^{-1} = (T\,M\,T^{-1})^{-1} = (T^{-1})^{-1}\,M^{-1}\,T^{-1}=T\,M^{-1}\,T^{-1}=(M^{-1})^\rho$.
\item $T\,M_1\,T^{-1}+T\,M_2\,T^{-1} = T\,(M_1+M_2)\,T^{-1}$ ~due to distributivity (similar for vectors).
\item Because of the associativity, it holds
      \begin{eqnarray}
      (T\,M_1\,T^{-1})\,(T\,M_2\,T^{-1}) &=& T\,M_1\,(T^{-1}\,T)\,M_2\,T^{-1} = T\,(M_1\,M_2)\,T^{-1} \mbox{~and} \nonumber\\ 
      (T\,M\,T^{-1})\,(T\,v)             &=& T\,M\,(T^{-1}\,T)\,v = T\,(M\,v). \nonumber
      \end{eqnarray} 
\item The proof is similar.
\end{enumerate}
\end{proof}

Other properties of reordering which are not directly related to the content of the paper are summarized in Appendix \ref{appreorder}.
For examples, see section~\ref{exreorder}.
\newpage

\subsection{State symmetries}

For the identity permutation $\rho=id$, the re-ordering leads naturally to equality: $v^{id}=v$ and $M^{id}=M$. 
If there are other permutations for which these equations hold, we call them ''proper state symmetries''.

\begin{definition}[State symmetries of vectors and matrices] ~
\begin{enumerate}
\item A permutation $\rho$ is called a \underline{state symmetry} of the vector $v$, iff $v=v^\rho$.
\item A permutation~$\rho$ is called a \underline{state symmetry} of the matrix $M$, iff $M=M^\rho$.
\end{enumerate}
If $\rho\ne id$ then the state symmetry is called \underline{proper}.
If a vector or matrix has a~proper state symmetry, it is called \underline{state symmetric}.
\end{definition}
With other words: state symmetric vectors $v$ and matrices $M$ fulfil the condition $v=T\,v$ or
$M=T\,M\,T^{-1}$, respectively, where $T$ is a permutation matrix.

\smallskip

As a first corollary, we characterise state symmetric objects:

\begin{corollary}[Equal elements in state symmetric objects] ~ \label{firstcor}
\begin{enumerate}
\item For vectors, a necessary and sufficient condition for having a proper state symmetry is the equality 
of at least two elements.
\item If a matrix $M=(m_{i\,j})$ has a state symmetry $\rho$, 
then the vector $v=(m_{i\,i})$ of the diagonal elements of $M$ has this state symmetry.
\end{enumerate}
\end{corollary}
\begin{proof} 
The requirement of having equal elements is straightforward. 
(Note for matrices, that the set of diagonal elements of $M$ and $M^\rho$ are always equal.)

On the other hand, if two elements $v_k$ and $v_l$ (with $k \neq l$) are equal, 
then the permutation \mbox{$\rho=(k\;l)(1)(2)\ldots(N)$} is a proper state state symmetry of $v$ so that the sufficience is shown.

For matrices, a counterexample for non-sufficience is the case, where $m_{k\,k}=m_{l\,l}$ for some $k\neq l$, 
but all other elements of $M$ are different from each other.
\end{proof}

For matrices, it is necessary to have at least two equal elements in the diagonal for 
beeing state symmetric. Moreover, all elements in the corresponding rows and columns must be equal.

\medskip

It can be expected that the investigation of state symmetries in vectors is much simpler than in matrices.
We will see this again in most of the theorems in the following sections.
For examples, see section~\ref{exreorder}.

\subsection{The groups of state symmetries}
This section shows properties of state symmetries for given vectors or given matrices.
We use the fact that state symmetries are permutations which may form groups.
The symmetric group $\SyG_N$ will be the basis to show subgroup properties of state symmetries.

\medskip

From the definition of state symmetry, we can directly conclude:

\begin{theorem}[Groups of state symmetries of a vector or a matrix] ~ \label{groupstate}
\begin{enumerate}
\item The collection of all state symmetries of a vector $v$ is a group with respect to composition 
(denoted: $\SSy(v):=\{\rho\in\SyG_N | \,v^\rho=v\}$).
\item The collection of all state symmetries of a matrix $M$ is a group with respect to composition 
(denoted: $\SSy(M):=\{\rho\in\SyG_N | \,M^\rho=M\}$).
\end{enumerate}
\end{theorem}
\begin{proof} 
$\SSy(v)$ and $\SSy(M)$ are isotropy subgroups of $v$ respective $M$.
\end{proof}

The statement for state symmetric matrices was first mentioned in (\cite{Behrends:00}, p.~17).

\medskip

We add a simple corollary on state symmetries of several vectors and matrices:

\begin{corollary}[Groups of state symmetries of several objects] \label{compgroup}~
\begin{enumerate}
\item The set of state symmetries of a set of column vectors $L_v=\{v^{(1)},v^{(2)}, \ldots\}$ is a group,
which is calculated by $\SSy(L_v)=\cap_{v^{(i)}\in L_v} \SSy(v^{(i)})$.
\item The set of state symmetries of a set of square matrices \mbox{$L_M=\{M^{(1)},M^{(2)}, \ldots\}$} 
is a group, which is calculated by \mbox{$\SSy(L_M)=\cap_{M^{(i)}\in L_M} \SSy(M^{(i)})$}.
\end{enumerate}
\end{corollary}

In addition, the application of a state symmetry to its permutation matix leads to:
\begin{corollary} 
Let $\rho$ be a permutation and $T$ the corresponding permutation matrix. Then $\rho$ is a state symmetry of $T$.
\end{corollary}
\begin{proof} 
$T=T\,(T\,T^{-1})=T^\rho$~. 
\end{proof} 

\pagebreak[3]

Example~5.2.\ref{expure} shows, that not for every group~$G$ there exist
vectors~$v$ or matrices~$M$ so that \mbox{$G=\SSy(v)$} or $G=\SSy(M)$. 
The next definition shall specify the groups for which this condition is fulfilled.

\begin{definition}\label{puress}
Let $\GG$ be a group of permutations of a state space $C$ with $|C|=N$.
\begin{enumerate}
\item If there exist a $N$-dimensional vector $v$ so that $\GG=\SSy(v)$, 
then $\GG$ is called a \underline{group of pure vector state symmetries}. 
\item If there exist a $N\times N$ matrix $M$ so that $\GG=\SSy(M)$, 
then $\GG$ is called a \underline{group of pure matrix state symmetries}. 
\end{enumerate}
\end{definition}

Because of Corollary~\ref{firstcor}.2, all groups of pure vector state symmetries
are groups of pure matrix state symmetries.

Groups of pure vector state symmetries can be easely characterised, using the symmetric groups $\SyG(C_k)$
of the elements of a partition of $C$.
\begin{lemma}\label{charvecgroup}
A group $\GG$ is a group of pure vector state symmetries on state space~$C$, iff it has the structure
\begin{equation}
\GG=\SyG(C_1)\times \SyG(C_2)\times \ldots \times \SyG(C_j)
\mbox{~,~ where ~} {\cal C}=\{C_i\}_i \mbox{~is a partition of $C$}. \label{form1}
\end{equation}
\end{lemma}
\begin{proof}
If $\GG$ is a group of pure vector state symmetries, then there is a vector $v$ so that 
$\GG=\SSy(v)$. Let $Q$ be the set of distinct elements of $v=(v_j)$. Then $C_i=\{j|q_i=v_j\}$.\\
On the other hand, if a group has the structure (\ref{form1}), then let 
\mbox{$v=(v_j=\sum i\,{\mathbf 1}_{j\in C_i}$)}. Then, $\GG=\SSy(v)$.
\end{proof}

For groups of pure matrix state symmetries, no such characterisation has been found yet.
\newpage
\section{The structure of state symmetric objects}

In this section, we investigate the opposite view. Here, the structure of vectors and matrices 
for given state symmetries shall be examined to show,
that they form subspaces of the vector spaces~$\Rspace$ and $\Mspace$
and subgroups of the linear group $GL_N({\mathbb R})$.

\subsection{Sums of of state symmetric objects}

This section derives properties for sums of state symmetric matrices and vectors. 

\begin{theorem}[Vector spaces on a single permutation] ~ \label{vecspacesym}
Let $\rho$ be a permutation on the set~ $C=\{1,2,\ldots,N\}$.
\begin{enumerate}
\item Let $\maV_\rho:=\{v\in \Rspace | \,v^\rho=v\}$ 
be the subset of~ $\Rspace$ of all vectors for which $\rho$ is a state symmetry.
Then $\maV_\rho$ is a vector space.
\item Let $\maM_\rho:=\{M\in \Mspace | \,M^\rho=M\}$ 
be the set of all square matrices for which $\rho$ is a state symmetry.
Then $\maM_\rho$ forms a vector space.
\end{enumerate}
\end{theorem}
\begin{proof} The subspace properties of $\maV_\rho$ and $\maM_\rho$ 
follow directly from Lemma~\ref{lemreorder}.3. and \ref{lemreorder}.5.
\end{proof}

In the previous chapter it was shown that state symmetries form groups, and that the groups
of several state symmetric objects can be derived directly.

Analogously to corollary~\ref{compgroup}, a statement for state symmetric objects can be concluded:

\begin{corollary}[Vector spaces on state symmetric objects] ~ \label{corvecspace}\\
Let $P=\{\rho_1,\rho_2,\ldots\}$ a set of permutations of $C=\{1,\ldots,N\}$.
\begin{enumerate}
\item The set of vectors $\maV\!_P\subset\Rspace$ for which exactly all permutations of $P$ are state symmetries is a vector space,
 which is calculated as $\maV\!_P=\bigcap_{\rho_i\in P} V_{\rho_i}$ . 
\item The set of square matrices $M_P\subset\Mspace$ for which exactly all permutations of $P$ are state symmetries
is a vector space, which is calculated as $\maM_P=\bigcap\limits_{\rho_i\in P} \maM_{\rho_i}$ .\\
\end{enumerate}
\end{corollary}
The following lemma simplifies the structure of the resulting vector spaces.
\begin{lemma}[Single state symmetry for spaces of state symmetric vectors] ~ \label{corvecspace1}\\
For every set of permutatations $P=\{\rho_1,\rho_2,\ldots\}$ there is a $\hat{\rho}\in\SyG(C)$ 
so that $\maV_P=\maV_{\hat{\rho}}$.
\end{lemma}
\begin{proof}
Constructing $\hat{\rho}$ from $P$: Two elements $c_i$ and $c_j$ are in the same cycle of $\hat{\rho}$ 
iff there is a $k$ so that $c_i$ and $c_j$ are in the same cycle of $\rho_k\in P$.
\end{proof}
Note that $\hat{\rho}\in <\!P\!>$ is not necesserely fulfilled (see example 5.3.\ref{otherperm}).

The dimension of $\maV_\rho$ is given by:
\begin{corollary}[Dimension of the state symmetric vector space] \label{dimvec}
For a given permutation $\rho$ of a state space $C$, the dimension of $\maV_\rho$ 
is the number $r$ of disjoint cycles of~$\rho$. 
\end{corollary}
A canonical basis of  $\maV_\rho$ can be derived using the partition of lemma~\ref{charvecgroup}.

\medskip

It is obvious that for every group $\GG$ on $C$ it holds $\GG\subseteq \SSy(\maV\!\sGG)$ and
\mbox{$\GG\subseteq \SSy(\maM\!\sGG)$}.
For groups of pure state symmetries, strict equality can be concluded:

\begin{theorem}[Completeness of groups of pure state symmetries] ~ \\ \label{completeness}
Let $C$ be a state space and $\GG$ be a group acting on $C$.
Let $\maV\!\sGG$ be the vector space of column vectors which have all state symmetries $\rho\in\GG$.\\
\mbox{$\GG=\SSy(\maV\!\sGG)$}, iff $\GG$ is a group of pure vector state symmetries.
\end{theorem}
\begin{proof} ~
Let $\GG$ be a group of pure vector state symmetries. Then there is a vector \mbox{$w\in\maV\!\sGG$}
so that $\GG=\SSy(w)$. Let $D=\{d_1,\ldots,d_j\}$ be the set of different elements of $w$, 
then $\{(v^{(l)}_i)\}_{l=1,\ldots,j} $ with $v^{(l)}_i=\mbox{\bf 1}_{l=d_i}$
is a set of basis vectors of $\maV\!\sGG$, so that \mbox{$w=\sum_{l=1}^j d_i\,v^i_l$}.

From lemma~\ref{charvecgroup}, the structure of $\GG$ is given by
\mbox{$\GG=\SyG(C_1)\times \SyG(C_2)\times \ldots \times \SyG(C_j)$}. 
Hence, the partition ${\cal C}=\{C_j\}$ of $C$ leads to the structure of the vector space
\begin{equation} \nonumber 
\maV\!\sGG=\maV_{{\cal S}(C_1)}\oplus \maV_{{\cal S}(C_2)}\oplus \ldots \oplus \maV_{{\cal S}(C_j)}
\mbox{~,~ where ~} {\cal C}=\{C_j\} \mbox{~is a partition of~} C
\end{equation}
and $\maV_{{\cal S}(C_l)}$ are one-dimensional vector spaces which are orthogonal to each other.
\end{proof}
For any set of permutations $P$, the group $\GG=\SSy(\maV\!_P)$ is a group of pure vector 
state symmetries because of lemma~\ref{corvecspace1}. However, $<P>\subseteq\GG$; the equality
is not fulfilled in general (see examples~5.3.\ref{vecspaceex1} and~5.3.\ref{vecspaceex2}).
\medskip

Only the corresponding statement to corollary~\ref{dimvec} 
has been proved also for spaces of state symmetrics matrices:
\begin{lemma}[Dimension of the space of state symmetric matrices]
Let $\rho$ be a~permutation of $C$.
The dimension of $\maM_\rho$ is given by the sum of the greatest common divisor 
of the length of the cycles:
\begin{equation}\label{GCDeq}
\dim(\maM_\rho)=\sum_{0\,\le\,i,j\,\le\,n_r} {\mathit GCD}(n(c_i),n(c_j))
\end{equation}
where $n_r$ is the number of cycles of $\rho$, and $n(c_i)$ is the length of the cycle $c_i$.
\end{lemma}
\begin{proof} Without loss of generalisation, $\rho$ can be reordered so that 
subsequent elements of $C$ are within each cycle of $\rho$.
Then, a matrix $M$ for which $\rho$ is a state symmetry can be constructed by dividing it into 
rectangular parts. 
In each part, the degree of freedom is $GCD(n(c_i),n(c_j))$, because not only $\rho$ but also 
powers of $\rho$ must be considered. The elements of each rectangle are independet 
from the elements of the other rectangles, so that they can be summed summed up directly,
leading to equation~(\ref{GCDeq}).
\end{proof}
For example, let  $\rho=(1\,2)(3\,4\,5)$, then
$M=\begin{pmatrix}  . & . & | & . & . & .  \\ 
                    . & . & | & . & . & .  \\ 
                    - & - & + & - & - & -  \\
                    . & . & | & . & . & .  \\ 
                    . & . & | & . & . & .  \\ 
                    . & . & | & . & . & .  \end{pmatrix}$, \\
which leads to
$M=\begin{pmatrix}  1 & 2 & 3 & 3 & 3  \\ 
                    2 & 1 & 3 & 3 & 3  \\ 
                    4 & 4 & 5 & 6 & 7  \\ 
                    4 & 4 & 7 & 5 & 6  \\ 
                    4 & 4 & 6 & 7 & 5  \end{pmatrix}$, 
so that $\dim(\maM_\rho)= 2+1+1+3= 7$ .

~
\medskip

The following conjecture (similar to theorem~\ref{completeness}) has not yet been proved.
\begin{hypothesis}[Completeness of groups of pure state symmetries for matrices]
Let $\maM\!\sGG$ be the vector space of matrices of a fixed dimension $N$ 
which have all state symmetries $\rho\in\GG$.
\mbox{$\GG=\SSy(\maM\!\sGG)$} iff $\GG$ is a group of pure matrix state symmetries.
\end{hypothesis}

For matrices, there has no hypothesis been found yet which 
determines the dimension of the space $\maM_P$ of the 
set~\mbox{$P=\{\rho_1,\rho_2\,\ldots\}$}
similar to lemma~\ref{corvecspace1}, where $\maM_P$ consists of all matrices 
which are state symmetric with regard to all permutations of $P$
(see example~5.3.\ref{corvecspace2}).\\ 
However, it is possible to derive the dimension algorithmically for any given set~$P$.

\medskip

In Appendix~\ref{basisvecex}, a canonical basis for all vector spaces  of sets $|C|\leq 4$ are presented.

\newpage

\subsection{Products of state symmetric objects}
\subsubsection{Groups of state symmetric matrices}

In this section, products of state symmetric objects will be investigated.

\begin{theorem}[Groups of state symmetric matrices] \label{symprodmat}
Let $C=\{1,2,\ldots,N\}$ and let $\rho\in\SyG(C)$ be a permutation of $C$.
Let~ ${\maM'_\rho} = \maM_\rho \cap GL_N({\mathbb R})$ be the set of all invertible square matrices 
for which $\rho$ is a state symmetry.
Then~ \mbox{${\maM'_\rho}\,$} forms a~group with respect to multiplication.
\end{theorem}
\begin{proof} The subgroup properties of~ ${\maM'_\rho} \subset GL_N({\mathbb R})$ can be shown
using Lemma~\ref{lemreorder}.2 and \ref{lemreorder}.4.
\end{proof}

Note also that if $M$ is state symmetric in respect to $\rho$, than also the power $M^k$ (for any
positive integer $k$). (If $M$ is invertible, then $k$ might indeed be a non-positive integer as well.)
\bigskip

Because $\maM'_\rho$ is acting on $\maV_\rho$ and $\maM_\rho$, we can also find statements for products
with state symmetric matrices and their factors. Invertibility is not in all cases required.

\begin{corollary}\label{symfac} 
Let $M_1$, $M_2$ and $M$ be $N\times N$ matrices and $v$ and $w$ $N$-dimensional column vectors,
with $M=M_1\,M_2$ and $w=M\,v$.
\begin{enumerate}
\item If $\rho$ is a state symmetry of the invertible matrix $M'$ 
      then $\rho$ is a state symmetry of $M_1$ iff it is a state symmetry of $M_2$.
\item If $\rho$ is a state symmetry of the invertible matrix $M_1$, 
      then $\rho$ is a state symmetry of $M_2$ iff it is a state symmetry of $M'$.
\item If $\rho$ is a state symmetry of the matrix $M$ and  vector $v$, 
      then it is a state symmetry of the product $w=M\,v$.
\item If $\rho$ is a state symmetry of the invertible matrix $M$, 
      then $\rho$ is a state symmetry of $w$ iff it is a state symmetry of $v$.
\end{enumerate}
\end{corollary}
The proof of 3. can be done using Lemma~\ref{lemreorder}.4. 
The other statements follow directly from the action of the group $\maM'_\rho$ 
on the vector spaces $\maV_\rho$ and $\maM_\rho$.\\

Corollary~\ref{symfac} is especially useful for the application of state symmetries to 
Markov chains (see Section~\ref{secMK} and \cite{Ring:02}).


\subsubsection{Sequences with state symmetric vectors}

Simple examples show that the existance of a state symmetry of the vectors $v$ and~$v'$ 
(with $v'=M\,v$) is not sufficient for the existance of a state symmetry of the matrix $M$ 
(see Example~5.4.\ref{exstatfac}). On the other hand, if $\rho$ is a state symmetry of both $v\neq 0$ and~$v'$,
then there is a matrix $M$ with $v'=M\,v$ and $\rho$ is a state symmetry of $M$ (for example, if
all elements of $v$ are different from zero, then a~diagonal matrix with diagonal elements 
$m_{ii}=v'_i/v_i$ has this state symmetry).

\medskip

In the following section, it shall be investigated which conditions are
necessary for the deduction of state symmetry of a matrix which generates a vector~sequence
$v^{(k)}=M^k\,v^{(0)}$.
From corollary~\ref{symfac}.3 we can conclude, that if $M$ is invertible and $M$~and $v^{(0)}$ 
have a common state symmetry, then all vectors of the gene\-rated sequence 
are state symmetric. 
Here, more precise prerequirements shall be presented which characterise 
generators of sequences of state symmetric vectors. First we define 
the term matrix generated vector sequences and derive some more general 
properties, which the will lead to the desired relation to state symmetries.
\smallskip

\begin{definition}[Matrix generated vector sequence]
Let $u=(v^{(k)})_{k\in {\mathbb N}}$ be a~\mbox{countable} infinite set of $N$-dimensional vectors. 
Let \mbox{$u(i)=(v^{(k)})_{k=0,\ldots,i}$} be the finite set of the first $i+1$ elements of $u$.
Moreover, let 
\begin{equation}
\mathbb{W}_{u(i)}:=\{M'|v^{(k+1)}=M'\,v^{(k)} \mbox{~for~all~} k=0,\ldots,i-1\}
\end{equation}
be the set of matrices which may generate $u(i)$. \\
Then, $u$ is a \underline{matrix generated vector sequence}, iff $\mathbb{W}_{u}\neq \emptyset$.
For any matrix \mbox{$M\in\mathbb{W}_{u}$}, $u$~may be called a $M$-generated vector sequence. 
\end{definition}
Hence, Markov chains are special $P$-generated vector sequences, where $P$ is a~probability matrix 
and $v^{(0)}$ is a probability vector.

\medskip

It is obvious that for every $i$ it holds that $\mathbb{W}_{u(i)}\supseteq\mathbb{W}_{u(i+1)}$, 
so that the sets $\mathbb{W}_{u(i)}$ form a non-increasing chain. Note that for a matrix generated 
vector sequence, if $v^{(k)}=v^{(k+1)}$ for one $k$, then $v^{(k)}=v^{(k+l)}$ for any $l\ge 0$. 
In particularely, if $v^{(k)}=0$ for one $k$ then $v^{(k+1)}=v^{(k)}=0$ and $v^{(k+l)}=0$ for any $l\ge 0$.

\begin{proposition}\label{proposelin}
Let $u(N)$ be a sequence of $N$-dimensional vectors. 
Iff the vectors of $u(N-1)$ are \mbox{linearely} independent, then
there is \mbox{exactly} one matrix $M$ so that $u(N)$ is M-generated.
\end{proposition}
\begin{proof}
Let $\maV_1=\maV_2=R^N$. The vectors $(v^{(i)})_{i=0,\ldots,N-1}$ are a basis of $\maV_1$. \\
Moreover, the vectors $(v^{(i)})_{i=1,\ldots,N}$ are all in $\maV_2$. Hence, there is exactly 
one map $\maV_1\to\maV_2$, which can be which can be calculated in terms of a matrix $M$ as 
$M=B\,S^{-1}$, where $S$ is the matrix of column vectors $(v^{(i)})_{i=0,\ldots,N-1}$ and 
$B$ is $(v^{(i)})_{i=1,\ldots,N}$. 
If the vectors of $u(N-1)$ are linearely dependent, then $S$ is not invertible.
Hence, for every sequence of $N$-dimensional vectors, \mbox{$|\mathbb{W}_{u(N-1)}|\neq 1$}.
\end{proof}

\begin{proposition}\label{corlinnumb}
If $u=(v^{(i)})_i$ is matrix generated vector sequence on $\Rspace$,
then the vectors of $u(N)$ are linearely dependent. 
\end{proposition}

These propositions imply that a matrix generated vector sequence $u$ 
is determined by its first $N+1$ elements $u(N)$. In the following, is it shown that
this number can be reduced in special cases.

\begin{lemma}\label{hyplin}
Let $u=(v^{(i)})_i$ be a matrix generated vector sequence. 
Let $i^*$ be defined as the number, for which the vectors of $u(i^*-1)$ are linearely independent,
but the vectors of $u(i^*)$ are linearely dependent.
Then, for every two matrices \mbox{$M_1,M_2\in\mathbb{W}_{u(i^*)}$} it holds that 
$M_1 \,v^{(i^*)}=M_2 \,v^{(i^*)}$.
\end{lemma}
\begin{proof}
Because of the linear dependence, there exists numbers $a_0,\ldots,a_{i^*}$ which not all are equal 
to zero so that $0=\sum_{i=0}^{i^*} a_i\,v^{(i)}$. For the case $v^{(i^*)}=0$, the claim is trivial.
Otherwise, $a_{i^*}\neq 0$. Without loss of generalisation, we set $a_{i^*}=1$.

Let $M$ be a matrix so that $v^{(i+1)}=M\,v^{(i)}$ for $i=0,\ldots, i^*$. 
Then
\begin{eqnarray}
v^{(i^*+1)} &=& M\,v^{(i^*)} \nonumber \\ \nonumber
            &=& M\,\sum_{i=0}^{i^*-1} a_i\,v^{(i)} \\ \nonumber
            &=& \sum_{i=0}^{i^*-1} a_i\, M\,v^{(i)} \\ \nonumber
            &=& \sum_{i=0}^{i^*-1} a_i\, v^{(i+1)} 
\end{eqnarray}
so that continuation of the vector sequence $u$ only depends on the numbers $a_0,\ldots,a_{i^*-1}$ and
not on a specific generator $M$ of $u$.
\end{proof}
As a direct conclusion, we obtain:
\begin{corollary}\label{corsetchain}
Let $u=(v^{(i)})_i$ be a matrix generated vector sequence and let \mbox{$m\in{\mathbb N}$}. 
If the vectors of $u(m)$ are linearely independent, then 
$\mathbb{W}_{u(m+1)}\subsetneqq \mathbb{W}_{u(m)}$, otherwise
$\mathbb{W}_{u(m)}= \mathbb{W}_{u(m+1)}=\mathbb{W}_{u}$.
\end{corollary}
Lemma~\ref{hyplin} and corollary~\ref{corsetchain} imply that the chain of 
sets $(\mathbb{W}_{u(m)})_m$ is strictly decrea\-sing for $m\le i^*$, while for $m>i^*$ all elements
are equal to the limit $\mathbb{W}_u$. 

As a consequence, a matrix generated vector sequence is fixed by the first~$i^*\le N$ vectors, 
if $i^*$ is defined as above. There might be several matrices $M\in\mathbb{W}_u$ which generate $u$,
but the general relationship between the vetcors of $u$ is determined by the linear dependence 
asshown in Lemma~~\ref{hyplin}.

This gives raise to the following definition.
\begin{definition}
For a matrix generated vector sequence with $i^*$ defined as above, 
the numbers $a_0,\ldots,a_{m-1}$ are called the \underline{linear coefficients} of $u$.
(In the case $v^{(i^*)}=0$, the linear coefficients are defined as zero.)
\end{definition}

\medskip

Note that in the case $i^*=N$, the linear coefficients are determined by the following corollary:
\begin{corollary}\label{corrcharpoly}
Let $u=(v^{(i)})_i$ be a matrix generated sequence of $N$-dimensional vectors.
If $i^*=N$, then the coefficients of the characteristic polynomial 
\mbox{$P_M(t)=(M-t\,I)=c_0+c_1\,t+\ldots+c_N\,t^N$} coincide with the linear 
coefficients $a_i$ of~$u$ (with exception to the sign of $c_N=(-1)^N$, which is multiplied to
the other coeffients: $a_i=c_N\,c_i$ for $i=0,\ldots,N-1$).
In particularely, $a_0=|M|$ and \mbox{$a_{N-1}=\mbox{trace}(M)$}.
\end{corollary}
\begin{proof}
For every $M$-generated sequence $(v_i)_i$ it holds:
\begin{eqnarray}
\sum_{i=1}^N c_i\,v_i &=& \sum_{i=1}^N c_i\,(M^i\,v_0)              \nonumber\\
                      &=& \left(\sum_{i=1}^N c_i\,M^i\right) \,v_0  \nonumber\\
                      &=& P_M(M) \,v_0                              \nonumber\\
                      &=& (M-M\,I) \,v_0                            \nonumber\\
                      &=&  0                                        \nonumber
\end{eqnarray}
and the claim is shown because $c_N=(-1)^N$ and $i^*=N$.
\end{proof}
The derivation above is also true for $i^*<N$, but in this case the coefficients of the
characteriatic polynomial~$c_i$ are not the linear coefficients of~$u$.
\bigskip

The final theorem summarizes these relations for matrix generated vector sequences 
and state symmetric matrices. We define the set of generators with state symmetry~$\rho$
as ${\mathbb W}^\rho_u = {\mathbb W}_u \cap \maM_\rho$ .

\begin{theorem}\label{maintheor}
Let $\maV=\Rspace$ a vector space of dimension $N>1$ and let 
\mbox{$u=(v^{(k)})_{k\in {\mathbb N}}$} be a matrix generated vector sequence. 
Let $\rho$ be a permutation of $C=\{1,\ldots,N\}$ which consists of $r$ cycles.
If $\rho$ is a common proper state symmetry of all vectors of $u(r)$, then the vectors 
of $u(r)$ are linearely dependent, and $\rho$ is a state symmetry of all vectors of $u$.
\end{theorem}
\begin{proof}
Corollary~\ref{dimvec} shows that $dim(V_\rho)=r$. Because $u(r)\in V_\rho$, 
proposition~\ref{corlinnumb} proofs the theorem.
\end{proof}
The main statement of this theorem is that the number of required vectors to determine $u$
can be reduced from from $N$ to $r$, based on the state symmetric property.

\medskip

\begin{corollary}[State symmetric generator]
Let $\maV=\Rspace$ a vector space of dimension $N>1$ and let \mbox{$u=(v^{(k)})_{k\in {\mathbb N}}$} 
be a matrix generated vector sequence. 
If $\rho$ is a common proper state symmetry of all vectors of $u$, then $|\mathbb{W}_u|>1$
and $|\mathbb{W}^\rho_u|\ge 1$.
\end{corollary}

Note that the reversion of this corollary above is not true in general 
(see example~5.5.\ref{reversalfalse}).
\pagebreak  

\section{Applications}

In this section, we would like to introduce briefly the intentions for applications of state 
symmetries in other mathematical fields, namely in Marov chains and in graph theory. 
For details, we refer to future papers of the autors which are in preparation.

\subsection{Markov chains} \label{secMK}

As mentioned in the introduction, we see the most important applications within the area 
of Markov chains.
Special structures of the transition matrix may be used to reduce the state space 
and which might decrease the time and effort for (numerical) calculations. 
This methodology is called aggregation of Markov chains or ''lumping''; 
it has been developed since the publication of \cite{Rosenblatt:57}. 
A large number of results related to matrices and vectors had been summarized 
in \cite{Kemeny:60}. However, first \cite{Benes:78} used group theory in a similar manner
like our approach, but only applied them to the limit distribution of the Markow chain.
\medskip

The definition of ''strong lumpability'' and ''weak lumpability'' follows the monograph of 
\cite{Kemeny:60}, and the exact lumpability has been defined by \cite{Schweitzer:84}.
\begin{definition}
Let $X=(P,v_0)$ be a Markov chain on state space $C$ with transition matrix $P$ and start vector $\pi$.
Let ${\cal C}=\{C_1,C_2,\ldots,C_r\}$ be a partition of $C$ and~
$f_A:\; f_A(c)=\sum_i i\,{\mathbf 1}_{c\in C_i}$ be a surjective map from $C$ to $\tilde{C}$.
\begin{itemize}
\item
The Markov chain $X$ is strongly lumpable with respect to the partition ${\cal C}$
if for every starting vector $\pi$ the lumped process $f_A(X)$ is a Markov chain and the 
transition probabilities do not depend on the choice of $\pi$.
\item
The Markov chain $X$ is weakly lumpable with respect to the partition ${\cal C}$
if there exists a starting vector $\pi$ so that the lumped process $f_A(X)$ is a~Markov chain.
\item 
The Markov chain $X$ is exactly lumpable with respect to the partition ${\cal C}$
if~for all states $s_k\in C_i$ of the aggregated process $f_A(X)$ holds at any time
\[
\Prob[X_n=s_k]=\frac{\Prob[\tilde{X}_n=r]}{|C_r|} \,.
\]
\end{itemize}
\end{definition}
If a Markov chain is exactly lumpable or strongly lumpable then it is weakly lumpable.

An example for the application of lumpability is the theorem of \cite{Barr:77}. 
It states, that the Eigenvalues of transition matrices on the aggregated state space are all found 
in the transition matrix on the original state space, if the Markov chain is strongly lumpable.
\medskip

Similar to the definitions above, one may define ''symmetric lumpability'' \cite{Ring:02}. 
This is based on the following theorems:
\begin{theorem}[Strong lumpability of state symmetric transition matrices] ~\\
Let $X=(P,v_0)$ be a Markov chain where $P$ is a transition matrix with state symmetry $\rho$.\\
Then, $X$ is strongly lumpable with respect to the associated partition of~$\rho$.
\end{theorem}
\begin{theorem}[Exact lumpability of state symmetric Markov chains] ~\\
Let $X=(P,v_0)$ be a Markov chain where $P$ and $v_0$ have a common state symmetry $\rho$.
Then, $X$ is exactly lumpable with respect to the associated partition of~$\rho$.
\end{theorem}

The symmetric lumpability of Markov chains is a specification of existing aggregations 
methods which largely preserves information of the original system: The time course of the original 
chain can be reconstructed from the lumped chain, and the 1-distance and relative entropy divergence 
of the chain from its limit distribution are invariant under the lumping lumping process \cite{Ring:02}.

\medskip

\begin{figure}[ht]
\centering
\includegraphics[height=4cm,clip]{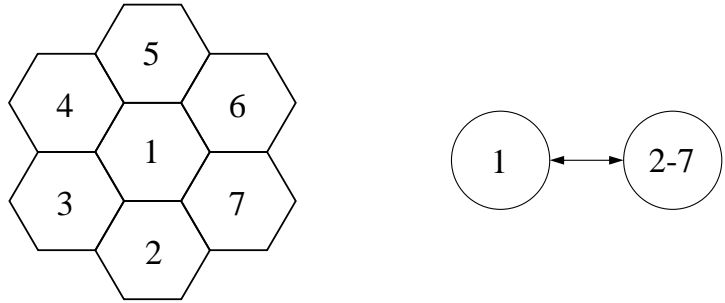}
\caption[The structure of the Markov chains]{The structure of the Markov chain $X$ 
(on the left) and the aggregated Markov chain $X'$ (on the right). For $X$, all transition 
probabilities between neighbouring cells are equal.\label{sixtwocells}}
\end{figure}

{\bf Example} ~
As an example, let $X=(P,v)$ be a Markov chain where 
\[P=\begin{pmatrix} 1-6\,p &   p    &   p    &   p    &   p    &   p    &   p \\ 
                      p    & 1-3\,p &   p    &   0    &   0    &   0    &   p \\
                      p    &   p    & 1-3\,p &   p    &   0    &   0    &   0 \\
                      p    &   0    &   p    & 1-3\,p &   p    &   0    &   0 \\
                      p    &   0    &   0    &   p    & 1-3\,p &   p    &   0 \\
                      p    &   0    &   0    &   0    &   p    & 1-3\,p &   p \\
                      p    &   p    &   0    &   0    &   0    &   p    & 1-3\,p \end{pmatrix}\]
and \mbox{$v=(1,0,0,0,0,0,0)\Transp$} (see figure~\ref{sixtwocells}).

The group $\cal G$ of state symmetries of $P$ is generated by $\rho_1=(1)(2\,3\,4\,5\,6\,7)$ and
$\rho_2=(1)(2\,3)(4\,7)(5\,6)$, so that it has the structure of the dihedral group ${\cal D}_6$.
All elements of $\cal G$ are also state symmetries of $v$.

$X$ is strongly lumpable with respect to the partition \mbox{${\cal C}=\{\{1\},\{2,3,4,5,6,7\}\}$}, 
and the aggregated process $X'=(P',v')$ is defined by the transition matrix 
\mbox{$P'=\begin{pmatrix}  1-6\,p &  p  \\ 
                      p     & 1-p \end{pmatrix}$}
and the starting vector $v'=(1,0)$.

Moreover, $X$ is exactly lumpable with respect to $\cal C$. 
Because of theorem~\ref{maintheor}, $X$ is state symmetric at any time.

\subsection{State symmetric graphs}

\paragraph{Weighted graphs}
As transition matrices of Markov chains can be interpreted as weighted graphs, 
the idea of state symmetries and agrregation can be transferred to graph theory.

Let $G=<V,E>$ be a connected, directed, weighted graph with the set of vertices $V$ and set of edges $E$.
$E$ can be described be the adjacency matrix $M_E$, and $V$ is the state space.
A state symmetry $\rho$ of the graph $G$ is a permutation of $V$ so that $E^\rho=E$.
The ''state symmetric lumping of a graph'' can be defined using the adjacency matrix and 
the strong lumpability aggregation from the previous section as shown in figure~\ref{sixtwocells}.

\paragraph{Unweighted graphs}

For unweighted (undirected) graphs, the adjacency matrix contains only the elements 0 and 1.
The investigation of unweighted graphs provides statements on the maximally possible
state symmetries of a given graph structure and, therefore, of the maximal aggregation of graphs.

\begin{figure}[H]
\centering
\includegraphics[height=40mm,clip]{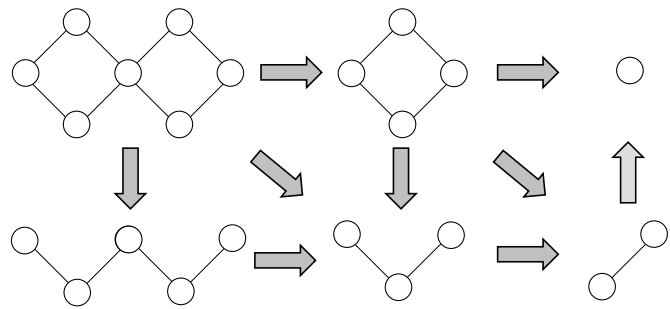}
\caption[State symmetric aggregation of graphs]{All possible state symmetric aggregations of a 
special graph with seven vertices. \label{stategraph}}
\end{figure}

There are new questions which can be raised for unweighted state symmetric graphs. For example, 
figure~\ref{stategraph} shows, that all of the five presented graphs can be aggregated to 
the simple 1-graph by a chain of state symmetric aggregations. An unsolved problem so far is 
for which structure the graph can be reduced to the 1-graph. It can be found that all graphs 
(without loops) with up to five vertices can be aggregated to the 1-graph. 

There is a graph with six vertices that cannot be aggregated to the 1-graph:
\begin{figure}[H]
\vspace*{1mm}
\centering
\includegraphics[height=20mm,clip]{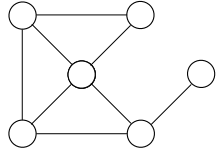}
\vspace*{3mm}
\caption[The graph with]{The graph with the smallest number of vertices that cannot be 
symmetrically aggregated to the \mbox{1-graph}.}
\end{figure}
If loops are not allowed, there exists one graph with seven vertices without state symmetries:
\begin{figure}[H]
\vspace*{1mm}
\centering
\includegraphics[height=15mm,clip]{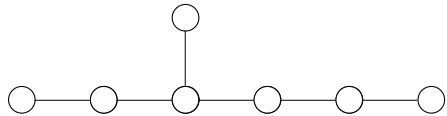}
\vspace*{3mm}
\caption[The smallest loop-free graph that cannot be symmetrically aggregated to the 1-graph]{The 
smallest loop-free graph that cannot be symmetrically aggregated to the \mbox{1-graph}.}
\end{figure}

\pagebreak

\section{Examples} \label{examples}

In this chapter, we would give examples (and counterexamples) to illustrate the statements of the 
preceeding chapters.

\subsection{Reordering and State Symmetries} \label{exreorder}
The first examples illutrates the concepts of reordering and state symmetries 
(see section~\ref{reordering}).
Let $C=\{1,2,3,4\}$ be the state space and let $\rho_1=(1\,3\,4)(2)$, $\rho_2=(1\,2\,3\,4)$ 
and \mbox{$\rho_3=(1)(2)(3\,4)$} be permutations of $C$.\\

\begin{asparaenum}[\bfseries{Example~\arabic{section}.\arabic{subsection}.}1]
\item (Permutations and their transformation matrices)
The permutation matrices of the permutations are 

     $T(\rho_1)=\begin{pmatrix}  0 & 0 & 1 & 0 \\ 
                                 0 & 1 & 0 & 0 \\
                                 0 & 0 & 0 & 1 \\
                                 1 & 0 & 0 & 0 \end{pmatrix}$, 
     $T(\rho_2)=\begin{pmatrix}  0 & 1 & 0 & 0 \\ 
                                 0 & 0 & 1 & 0 \\
                                 0 & 0 & 0 & 1 \\
                                 1 & 0 & 0 & 0 \end{pmatrix}$ and
     $T(\rho_3)=\begin{pmatrix}  1 & 0 & 0 & 0 \\ 
                                 0 & 1 & 0 & 0 \\
                                 0 & 0 & 0 & 1 \\
                                 0 & 0 & 1 & 0 \end{pmatrix}$.
\\~\\
\item (Reordering of matrices)

Let $M_1=\begin{pmatrix}         11 & 12 & 13 & 14 \\ 
                                 21 & 22 & 23 & 24 \\
                                 31 & 32 & 33 & 34 \\
                                 41 & 43 & 43 & 44 \end{pmatrix}$.
The resulting reordered matrices are
       $M_1^{\rho_1}=\begin{pmatrix} 33 & 32 & 34 & 31 \\ 
                                     23 & 22 & 24 & 21 \\
                                     43 & 42 & 44 & 41 \\
                                     13 & 12 & 14 & 11 \end{pmatrix}$, ~~
 \mbox{$M_1^{\rho_2}=\begin{pmatrix} 22 & 23 & 24 & 21 \\ 
                                     32 & 33 & 34 & 31 \\
                                     42 & 43 & 44 & 41 \\
                                     12 & 13 & 14 & 11 \end{pmatrix}$} and
  \mbox{$M_1^{\rho_3}=\begin{pmatrix} 11 & 12 & 14 & 13 \\ 
                                     21 & 22 & 24 & 23 \\
                                     41 & 42 & 44 & 43 \\
                                     31 & 32 & 34 & 33 \end{pmatrix}$}.
\\~\\
\item (Reordering and state symmetry of vectors) 

Let $v_1=(1,2,3,4)\Transp$ and $v_2=(1,6,1,1)\Transp$ be vectors on $C$.\\
Then, $v_1^{\rho_1}=(3, 2 , 4 , 1 )\Transp$;~~
      $v_2^{\rho_1}=(1, 6 , 1 , 1 )\Transp=v_2=v_2^{\rho_3}$;~~
      $v_2^{\rho_2}=(6, 1 , 1 , 1 )\Transp$, so that $\rho_1$ and $\rho_3$ are state symmetries of $v_2$.
      
\pagebreak[4]
\item (Reordering and state symmetry of matrices) 

Let $M_2=\begin{pmatrix}         1 & 2 & 3 & 4 \\ 
                                 5 & 6 & 5 & 5 \\
                                 4 & 2 & 1 & 3 \\
                                 3 & 2 & 4 & 1 \end{pmatrix}$ 
be a matrix on $C$, and let $v_2$ be as in the previous example.
(Note: $v_2$ is the vector of diagonal elements of $M_2$.) \\
Then, $M_2^{\rho_1}=M_2$;~~ 
  $M_2^{\rho_2}=\begin{pmatrix} 6 & 5 & 5 & 5 \\ 
                            2 & 1 & 3 & 4 \\
                            2 & 4 & 1 & 3 \\
                            2 & 3 & 4 & 1 \end{pmatrix}$ and
  $M_2^{\rho_3}=\begin{pmatrix} 1 & 2 & 4 & 3 \\ 
                            5 & 6 & 5 & 5 \\
                            3 & 2 & 1 & 4 \\
                            4 & 2 & 3 & 1 \end{pmatrix}$.\\
\end{asparaenum}

$\rho_1$ is a proper state symmetry of $M_2$ and $v_2$, while $\rho_2$ is not. 

$\rho_3$ is a state symmetry only of $v_2$, not of $M_2$.
It appears, that the degrees of freedom for the elements of state symmetric matrices and vectors 
is very limited (compare Correlary~\ref{firstcor}).

\subsection{Groups of state symmetries}

Let $M_2$, $v_2$, $\rho_1$ and $\rho_3$ be as in the previous subsection. 
\\
\begin{asparaenum}[\bfseries{Example~\arabic{section}.\arabic{subsection}.}1]
\item (Groups of state symmetries of a given vector/matrix)

The groups of state symmetries can be easily obtained using the state symmetries as generators:
\begin{eqnarray}
\SSy(M_2) &=& <\rho_1>=\{(1\,3\,4)(2),(1\,4\,3)(2),(1)(2)(3)(4)\}
\nonumber \\ \nonumber
\SSy(v_2) &=& <\rho_1,\rho_3>=\{(1\,3\,4)(2),(1\,4\,3)(2),(1)(2)(3)(4),\\ \nonumber
&& \qquad\qquad\qquad(1)(2)(3\,4),(1\,3)(2)(4),(1\,4)(2)(4)\} \\ \nonumber
&=& \SyG(\{1,3,4\})\times\{(2)\}
\end{eqnarray}

\item ~ Let
\[M_3=\begin{pmatrix}  1 & 2 & 3 & 4 & 4\\ 
                       3 & 1 & 2 & 4 & 4\\
                       2 & 3 & 1 & 4 & 4\\
                       5 & 5 & 5 & 6 & 7\\
                       5 & 5 & 5 & 7 & 6\end{pmatrix}\,.\]
The group $\SSy(M_3)$ contains of 6 elements and is generated by the permutation $<(1\,2\,3)(4\,5)>$:
\[ \SSy(M_3)=\SyG(\{1,2,3\}) \times \SyG(\{4,5\}) \,.\]
Let  $M_4=\begin{pmatrix}  1 & 2 & 2 & 3 & 3\\ 
                           2 & 1 & 2 & 3 & 3\\
                           2 & 2 & 1 & 3 & 3\\
                           4 & 4 & 4 & 5 & 6\\
                           4 & 4 & 4 & 6 & 5\end{pmatrix} $ .
The group $\SSy(M_4)$ contains of 12 elements and is generated by two permutations 
\mbox{$\SSy(M_4)=<(1)(2\,3)(4\,5),(1\,2\,3)(4)(5)> $} .

$M_4$ can be obtained from $M_3$ by decrementing all numbers larger than 2 \mbox{by 1}. 
Hence, the 3 becomes 2 and there are more symmetries possible, so that $\SSy(M_3)\subset\SSy(M_4)$.\\

\item\label{expure} (Groups of pure state symmetries) 
We show that there exist groups which are not ''groups of pure state symmetries''.

Let $\AlG_3$ be the alternating group of order 3, acting on $C=\{1,2,3\}$. 
Any vector $v=(a,b,c)\Transp$ for which $\SSy(v)=\AlG_3$ holds, must fulfil the condition \mbox{$a=b=c$}. 
In this case, also $\rho'=(1\,2)(3)\notin \AlG_3$ is a state symmetry of $v$, so that $\SSy(v)=\SyG_3$.
Hence, $\AlG_3$ is not a group of pure vector state symmetries. It is, however, a group of
pure matrix state symmetries, because for 
$M=\begin{pmatrix}  1 & 2 & 3 \\ 
                    3 & 1 & 2 \\
                    2 & 3 & 1 \end{pmatrix}$, $\AlG_3=\SSy(M)$.

Let $\AlG_4$ be the alternating group of order 4, acting on $C'=\{1,2,3,4\}$. 
Any matrix $M=(m_{i\,j})$ for which $\SSy(M)=\AlG_4$ holds, must fulfil the conditions 
\begin{eqnarray}
m_{ii} &=& m_{jj} \mbox{~~for all~} i       \mbox{~and~} j \mbox{~~~and} \nonumber \\ \nonumber
m_{ij} &=& m_{kl} \mbox{~~for all~} i\neq j \mbox{~and~} k\neq l \,.
\end{eqnarray}
Hence, $M=\begin{pmatrix}  a & b & b & b \\ 
                           b & a & b & b \\
                           b & b & a & b \\
                           b & b & b & a \end{pmatrix}$ 
for $a\neq b$, so that also \mbox{$\rho''=(1\,2)(3)(4)\notin \AlG_4$} is a~state symmetry of $M$.
Hence, $\SSy(M)=\SyG_4$ and $\AlG_4$ is not a group of pure matrix state symmetries.
\end{asparaenum}

\subsection{Vector spaces}

The vector spaces of matrices and vectors are presented in Appendix~\ref{basisvecex}
for all groups of state symmetries which are subgroups of the $\SyG_4$.
Here we give examples for special cases of Theorem~\ref{completeness} and Lemma~\ref{corvecspace1}.
\\
\begin{asparaenum}[\bfseries{Example~\arabic{section}.\arabic{subsection}.}1] 
\item ~ \label{vecspaceex1} 
Let ${C}=\{1,2,3\}$ and $\rho=(1\,2\,3)$.
Then \mbox{$<\rho>=\CyG_3=\{(1\,2\,3),(1\,3\,2),id\}$},
while $\SSy(\maV_\rho)=\SyG_3$. 
(The underlying partition contains of only one set, so that also $\rho'=(1\,2)(3)$ is a state symmetry.)
\\
\item ~ \label{vecspaceex2} Let ${C}$ and $\rho$ as in the previous example. 
Then, $\SSy(\maM_\rho)=<\rho>$ because the structure of any $M\in\SSy(\maM_\rho)$
is given by $M=\begin{pmatrix} a & b & c \\ 
                               c & a & b \\
                               b & c & a \end{pmatrix}$,
and there are no other state symmetries for these matrices in general.
\\
\item \label{otherperm}
It shall be shown that the vector $\hat{\rho}$ from Lemma~\ref{corvecspace1} does not necesserely 
fulfill $\hat{\rho}\notin <\!P\!>$.
Let $P=\{(1\,2)(3\,4),(1\,3)(2\,4)\}$, then $<\!P\!>$ is a Klein-4-group. 
However, all vectors that have both state symmetries of $P$ must consist of only one element 
($v=(a,a,a,a)\Transp$). Hence, $\hat{\rho}$ must have exactly one orbit 
(e.g. $\hat{\rho}=(1\,2\,3\,4)$) -- while all elements of $<\!P\!>$ have at least two orbits.
\\
\item \label{corvecspace2}
For matrices, there is no direct analogue statement to lemma~\ref{corvecspace1}, because
for any $k$ there is a set of permutations $P$ so that at least $k$ permutations are required to 
form $\maM_P$: Let  $P=\{\rho_1,\rho_2,\ldots,\rho_k\}$ with $\rho_i=((2\,i-1)(2\,i))$,
then there is no set $P'$ with less than $k$ elements so that $\maM_P=\maM_{P'}$.
\end{asparaenum}

\subsection{State symmetric factors} \label{exampfacvec}

This section presents examples and counterexamples for the state symmetry of products of 
matrices with vectors.
\\
\begin{asparaenum}[\bfseries{Example~\arabic{section}.\arabic{subsection}.}1]
\item \label{exstatfac}
Let $v^{(0)}=(1,1,-2)\Transp$ and  
\mbox{$M=\begin{pmatrix}   3 &  1 &  2  \\ 
                           0 & -2 & -1  \\
                           1 &  0 &  0  \end{pmatrix}$}.
Then \mbox{$v^{(1)}=M\,v^{(0)}=(0,0,1)\Transp$} is state symmetric but 
\mbox{$v^{(2)}=M^2\,v^{(0)}=(2,-1,0)\Transp$} is not.
$M$ is not state symmetric, so that corollary~\ref{symfac} cannot be applied.
\\
\item If we use instead a matrix with state symmetry $\rho=(1\,2)$, e.g.
\mbox{$M'=\begin{pmatrix}    1 & -5 & -2  \\
                            -5 &  1 & -2  \\
                   \frac{1}{2} &  \frac{1}{2} & 0  \end{pmatrix}$}, 
for which holds $v^{(1)}=M'\,v^{(0)}$ (similar as above),
then all elements of the sequence $(v^{(i)}=(M')^i\,v^{(0)})_i$ are state symmetric 
(e.g. \mbox{$v^{(2)}=(-2,-2,0)\Transp$} with $v^{(2)}=-2\,v^{(0)}-4\,v^{(1)}$).
The same vector sequence may, however, also be generated by non-state symmetric matrices, e.g.
$M''=\begin{pmatrix}          -5 & 1 & -2  \\
                             -2 & -2 & -2  \\
                              1 & 0  & 0  \end{pmatrix}$ . \medskip
\\
\item There are also examples of products of matrices: 
\[
M_1=\begin{pmatrix}   1 & 1   \\  2 & 2  \end{pmatrix} \qquad
M_2=\begin{pmatrix}   1 & -2  \\ -1 & 2  \end{pmatrix} \leadsto  
M=M_1\,M_2= \begin{pmatrix}   0 & 0  \\ 0 & 0  \end{pmatrix} 
\]
and while $M_1$ and $M_2$ are not state symmetric, the product $M$ is. 
The reason is that $M_1$ and $M_2$ are not invertible.
\\
\item
Let $M_3=\begin{pmatrix}   2 & 1  \\ 3 & -2  \end{pmatrix}$. 
$M_3$ is invertible, and its square is state symmetric although $M_3$ is not: 
$(M_3)^2=\begin{pmatrix}   7 & 0  \\ 0 & 7  \end{pmatrix}$. \\
Let $v_0=(a,a)\Transp$ be a state symmetrix vector, then $M_3\,v_0=(3\,a,a)\Transp$
is not state symmetric for $a\neq 0$, but $(M_3)^2\,v_0=(7\,a,7\,a)\Transp$ is.
\end{asparaenum}

\subsection{Sequences of state symmetric vectors} \label{exampseqvec}
\begin{asparaenum}[\bfseries{Example~\arabic{section}.\arabic{subsection}.}1]
\item Let $v^{(0)}=(0,1)\Transp$ and $M=\begin{pmatrix} 0 & 1 \\ 0 & 1 \end{pmatrix}$. 
Then, for $k>0$, all elements of the sequence $v^{(k)}=M^k\,v^{(0)}$ are constant $v^{(k)}=(1,1)$ 
and state symmetric. 
So a state symmetric sequence has been generated from a non state-symmetric vector, 
which is possible because $M$ is not state symmetric.
\\
\item Even simpler, if $M$ is the null matrix, then $w=M\,v=0$ is the (state symmetric) 
null vector for any $v$. In this example, $M$ is state symmetric but not invertible.
\\
\item 
Let $u=(v^{(i)})_i$ with $v^{(0)}=(1,1)\Transp$ and $v^{(1)}=(2,2)\Transp$. Then it follows: 
\begin{equation}
\mathbb{W}_{u(1)}=\left\{M\,\Big|\,M=\begin{pmatrix} a  &  2-a \\ 2-b & b \end{pmatrix}\right\}
\end{equation}
and for any vector $v=(c,c)\Transp$ and any $M\in\mathbb{W}_{u(1)}$, the product can be calculated as 
$M\,v=(2\,c, 2\,c)\Transp$. Hence, $\mathbb{W}_{u(2)}\neq\emptyset$ only if $v^{(2)}=(4,4)\Transp$, and in this case 
$\mathbb{W}_{u(1)}=\mathbb{W}_{u(2)}$. It follows that the sequence $u'=((2^i,2^i)\Transp)_i$ 
is generated by all elements of $\mathbb{W}_{u(1)}$. \\
Because the permutation $\rho=(1\,2)$ is a state symmetry of all elements of $u'$, 
there are state symmetric generating matrices $\bar{M}\in\mathbb{W}_{u'}$ with 
$\bar{M}= \begin{pmatrix} a  &  2-a \\ 2-a & a \end{pmatrix}$.
\\
\item \label{reversalfalse}
Let $u=(v^{(i)})_i$ with $v^{(0)}=(1,2,3)\Transp$, $v^{(1)}=(4,5,6)\Transp$ and \mbox{$v^{(2)}=(7,8,9)\Transp$}.
These vectors are linearely dependent with $v^{(2)}=-v^{(0)}+2\,v^{(1)}$, and hence, 
$v^{(3)}=-v^{(1)}+2\,v^{(2)}=(10,11,12)\Transp$.

The generating matrices of $u(2)$ have the structure:
\begin{equation}
\mathbb{W}_{u(2)}=\left\{M\,\Bigg|\,M=\begin{pmatrix} m_{11}    &  -1-2\,m_{11}  & 2+m_{11}   \\
											-1-m_{22}/2 &    m_{22}   & 2-m_{22}/2 \\
                                             m_{33}-4   & 5-2\,m_{33} &   m_{33}   \end{pmatrix}\right\}
\end{equation}
For any $M\in\mathbb{W}_{u(2)}$ it follows $M\,v^{(2)}=(10,11,12)\Transp$,
and $\mathbb{W}_{u(2)}=\mathbb{W}_{u(3)}$ (so that $|\mathbb{W}_{u(3)}|>1$).

There are two state symmetric matrices $\bar{M}\in\mathbb{W}_{u(3)}$,
\mbox{$\bar{M}_1= \begin{pmatrix} 0 & -1 & 2 \\
				                 -1 &  0 & 2 \\
						         -1 & -1 & 3 \end{pmatrix}$} with state symmetry $\rho_1=(1\,2)$ and
$\bar{M}_2= \begin{pmatrix} -1 & 1 & 1 \\
						    -2 & 2 & 1 \\
						    -2 & 1 & 2 \end{pmatrix}$ with state symmetry \mbox{$\rho_2=(2\,3)$}.
Hence, the reversal of theorem~\ref{maintheor} is not true, because in this example we have
$|\mathbb{W}_{u}|>1$ and $|\mathbb{W}^{\rho_1}_{u}|\ge 1$, but $\rho_1$ is no state symmetry of~$u$.
						    
Note in addition, that for $\rho_3=(1\,3)$ or  $\rho_4=(1\,2\,3)$, 
 $\mathbb{W}^{\rho_3}_{u(2)}=\mathbb{W}^{\rho_4}_{u(2)}=\emptyset$.
\\
\item Let $u=(v^{(i)})_i$ with $v^{(0)}=(1,1,0)\Transp$,  $v^{(1)}=(0,0,-1)\Transp$ and 
\mbox{$v^{(2)}=(-2,-2,-2)\Transp$} so that $v^{(2)}=-2\,v^{(0)}+2\,v^{(1)}$.
The generating matrices of $u(2)$ are:
\begin{equation}
\mathbb{W}_{u(2)}=\left\{M\,\Bigg|\,M=\begin{pmatrix} m_{11}  &  -m_{11} & 2 \\
											 -m_{22}  &   m_{22} & 2 \\
                                            -1-m_{32} &   m_{32} & 2  \end{pmatrix}\right\}
\end{equation}
There are infinitely many state symmetric matrices $\bar{M}\in\mathbb{W}_{u(2)}$, 
which have the structure 
$\bar{M}= \begin{pmatrix}  a      &       -a      & 2 \\
						  -a      &        a      & 2 \\
					 -\frac{1}{2} &  -\frac{1}{2} & 2 \end{pmatrix}$ with state symmetry $\rho=(1\,2)$.
\\
\item Let $u=(v^{(i)})_i$ be a matrix generated vector sequence with \mbox{$v^{(0)}=(1,1,1,2,2)\Transp$},
$v^{(1)}=(3,3,3,0,0)\Transp$ and $v^{(2)}=(4,4,4,2,2)\Transp$. 
These three vectors are linearely dependent and hence sufficient to fix $u$:
$v^{(2)}=v^{(0)}+v^{(1)}$ so that \mbox{$v^{(3)}=v^{(1)}+v^{(2)}=(7,7,7,2,2)\Transp$}.
\\
\item Let $u=(v^{(i)})_i$ be a matrix generated vector sequence with \mbox{$v^{(0)}=(1,2,3,4)\Transp$},
$v^{(1)}=(3,4,3,3)\Transp$, $v^{(2)}=(4,4,4,4)\Transp$, $v^{(3)}=(1,1,2,2)\Transp$ and $v^{(4)}=(2,2,3,3)\Transp$.
Then, $v^{(4)}=0\,v^{(0)}+0\,v^{(1)}+\frac{1}{4}v^{(2)}+v^{(3)}$ and 
\begin{equation}
\mathbb{W}_u=\mathbb{W}_{u(4)}=\left\{M=\begin{pmatrix} 
         -\frac{19}{4} &  \frac{13}{4} & \frac{23}{4}  &  -4 \\
         -\frac{19}{4} &  \frac{13}{4} & \frac{19}{4}  &  -3   \\
         -\frac{9}{2}  &  \frac{5}{2}  & \frac{15}{2}  &  -5  \\
         -\frac{9}{2}  &  \frac{5}{2}  & \frac{15}{2}  &  -5 \end{pmatrix}\right\} \,.
\end{equation}
$M$ is neither state symmetric nor invertible, and it is the only element of $\mathbb{W}_{u(4)}$.
However, for $i\ge 2$, all vectors of $v^{(i)}$ are state symmetric with $\rho_2=(1\,2)(3\,4)$
beeing  a state symmetry.
\medskip

If the first elements of $u$ are omitted, this generated sequence persists, 
but the set of generators becomes enlarged:
Let $u'=(v^{(i+1)})_i$ and $u''=(v^{(i+2)})_i$. Then, all vectors of $u'$ have the state symmetry
$\rho_1=(1)(2)(3\,4)$. The set of generators of $u'$ is given by
\begin{equation}
\mathbb{W}_{u'}=\mathbb{W}_{u'(3)}=\left\{\begin{pmatrix} 
         -\frac{19}{4} & \frac{13}{4} &  \frac{7}{4}-a  &   a  \\
         -\frac{19}{4} & \frac{13}{4} &  \frac{7}{4}-b  &   b  \\
         -\frac{9}{2}  & \frac{5}{2}  &  \frac{5}{2}-c  &   c  \\
         -\frac{9}{2}  & \frac{5}{2}  &  \frac{5}{2}-d  &   d \end{pmatrix}\right\}
\end{equation}
and the subset of state symmetric matrices $\mathbb{W}_{u'}^{\rho_1}$ is given by
\begin{equation}
\mathbb{W}_{u'}^{\rho_1}=\left\{\begin{pmatrix} 
         -\frac{19}{4}  &  \frac{13}{4} &  \frac{7}{8}    &  \frac{7}{8} \\
         -\frac{19}{4} &  \frac{13}{4} &  \frac{7}{8}    &  \frac{7}{8}  \\
         -\frac{9}{2}   &  \frac{5}{2}  &       a         &  \frac{5}{2}-a \\
         -\frac{9}{2}   &  \frac{5}{2}  &  \frac{5}{2}-a  &      a \end{pmatrix}\right\}\,.
\end{equation}

For $u''$ (with state symmetry $\rho_2=(1\,2)(3\,4)$), the generators are
\begin{equation}
\mathbb{W}_{u''}=\mathbb{W}_{u''(2)}=\left\{\begin{pmatrix} 
         -\frac{3}{2}-a &   a   & \frac{7}{4}-e  &  e     \\
         -\frac{3}{2}-a &   a   & \frac{7}{4}-f  &  f     \\
             -2-c       &   c   & \frac{5}{2}-g  &  g     \\
             -2-d       &   d   & \frac{5}{2}-h  &  h  \end{pmatrix}\right\}
\end{equation}
and the state symmetric generators are
\begin{equation}
\mathbb{W}_{u''}^{\rho_2}=\left\{\begin{pmatrix} 
              a         & -\frac{3}{2}-a &      c         &  \frac{7}{4}-c  \\
         -\frac{3}{2}-a &        a       &  \frac{7}{4}-c &     c           \\
              b         &      -2-b      &      d         &  \frac{5}{2}-d  \\
            -2-b        &        b       &  \frac{5}{2}-d &     d   \end{pmatrix}\right\} \,.
\end{equation}

~\\~
\item  {\bf $|\mathbb{W}_{u}| >1 $ does not imply that there is a $\rho$ 
so that $|\mathbb{W}^\rho_{u}|\ge 1$ .
}\\
Let $v^{(0)}=(1,1,0)\Transp$, $v^{(1)}=(1,0,-1)\Transp$ and $v^{(2)}=(0,1,1)\Transp$. 
These vectors are linearely dependent, and 
\begin{equation}
\mathbb{W}_{u}=\left\{\begin{pmatrix} 
         a  &  1-a &  a \\
        1+b & -1-b &  b \\
         0  &  -1  & -1 \end{pmatrix}\right\}\,.
\end{equation}
Hence $|\mathbb{W}_{u}|>1$, but there is no $\rho\neq id$ so that $\mathbb{W}^\rho_{u}\neq\emptyset$.
\\
\item 
{\bf It is not generally true that $\mathbb{W}^\rho_{u} \subsetneqq \mathbb{W}_{u}$ .}
\\
Let $v^{(0)}=(1,-1,2)\Transp$ and 
       $M=\begin{pmatrix} 
         1  & 2  & -1 \\
         2  & 1  & -1 \\
         0  & 0  & -2 \end{pmatrix}$ with state symmetry \mbox{$\rho=(1\,2)(3)$}.
Then, $v^{(1)}=(-3,-1,-4)\Transp$ and $v^{(2)}=(-1,-3,8)\Transp$.\\
Because $\det(v^{(0)},v^{(1)},v^{(2)})= -32$, these vectors are linearely independent, 
so that $\mathbb{W}_{u}=\mathbb{W}^\rho_{u}=\{M\}$. 
\end{asparaenum}

\subsection{Periodic vector sequences} \label{periodic}
Let $u=(v^{(i)})$ be a matrix generated vector sequence.
If $v^{(i)}\neq v^{(j)}$ for all $i\neq j$, then $u$ is called \underline{aperiodic}, 
otherwise it is periodic.
The definition of periodicity is meaningful, because if $v^{(i)}=v^{(j)}$ for $i\neq j$, 
then $v^{(i+n)}=v^{(j+n)}$ for every $n>0$, and $p=|i-j|$ is a period of $u$. 

\medskip

If $u$ is periodic with period $p$, then it is generated by a matrix $M$ for which holds
$M^p=I$.
\\
\begin{asparaenum}[\bfseries{Example~\arabic{section}.\arabic{subsection}.}1]
\item 
Let $v=(0,1)$ and $M=\begin{pmatrix} 0 & -1 \\ 1 & 1 \end{pmatrix}$. 
The sequence $v^{(k)}=v\,M^k$ is periodic with six different elements, because for every $l\in{\mathbb N}$
\begin{eqnarray}
v^{(6\,l)}   &=& (0,1)  \nonumber\\
v^{(6\,l+1)} &=& (1,1)  \nonumber\\
v^{(6\,l+2)} &=& (1,0)  \nonumber\\
v^{(6\,l+3)} &=& (0,-1) \nonumber\\
v^{(6\,l+4)} &=& (-1,-1)\nonumber\\\nonumber
v^{(6\,l+5)} &=& (-1,0)
\end{eqnarray}
Hence, $v^{(6\,l+1)}$ and $v^{(6\,l+4)}$ are state symmetric, while the other four elements are not.
The generating $M$ is invertible but not state symmetric.

\end{asparaenum}

\pagebreak  
\begin{appendix}
\section{Appendix}

\subsection{Additional properties of reodering}
\label{appreorder}
\begin{lemma}
Let $C$ be a finite state space, $\rho$ be a permutation of $C$ (with the \mbox{permutation} matrix $T$), 
$M_k=\left(m_{i\,j}\right)$ two $N\times N$ square matrices and $v_k$ two $N$-dimensional vectors ($k=1,\,2$). \\
Let $||.||_1$, $\!||.||_2$, $\!||.||_\infty$, $\!||.||_{\max}\,$ and $||.||_F$ be 
the absolute, euklidic, infinity, maximum and the Frobenius norm, defined as usual.
Then the following properties hold:
\begin{enumerate}
\item The determinants and the characteristic polynomials of $M$ and $M^\rho$ are equal, and
all Eigenvalues of $M$ are Eigenvalues of $M^\rho$ .
\item $v$ is an Eigenvector of $M$ iff $v^\rho$ is an Eigenvector of $M^\rho$ .
\item $||M||_1=||M^\rho||_1$ , ~ ~ ~ $||M||_2=||M^\rho||_2$ , ~ ~ ~ $||M||_\infty=||M^\rho||_\infty$ , ~ ~ ~
$||M||_{\max}=||M^\rho||_{\max}$ ~ and ~  $||M||_F=||M^\rho||_F$.
\item $||v||_1=||v^\rho||_1$ , ~ ~ ~ ~ $||v||_2=||v^\rho||_2$  ~ and ~ $||v||_\infty=||v^\rho||_\infty$ .  
\end{enumerate}
\end{lemma}
\begin{proof}
\begin{enumerate}
\item This is due to the similarity of $M$ and $M^\rho$.
\item If $v$ is an Eigenvector of $M$, then $\lambda\, v=M\,v$ and
  \begin{eqnarray}
    T\,\lambda\,v   &=& T\,M\,v \nonumber\\\nonumber
    \lambda\,T\,v   &=& T\,M\,(T^{-1}\,T)\,v = (T\,M\,T^{-1})\,(T\,v) \\\nonumber
    \lambda\,v^\rho &=& M^\rho\,v^\rho \,\,.
  \end{eqnarray}
\item For the euclidic norm $||.||_2$ and the Frobenius norm $||.||_F$, 
the statement is true because of the orthogonality of $T$. For the other norms we have:
  \begin{eqnarray}
  ||M||_1\,    &=& \max_{j\in C} \sum_{i=1}^N |m_{i\,j}| 
                = \max_{j\in C} \sum_{k=1}^N |m_{\rho(k)\,\rho(j)}|
                = ||M^\rho||_1 \nonumber 
\\
  ||M||_\infty &=& \max_{i\in C} \sum_{j=1}^N |m_{i\,j}| 
                = \max_{i\in C} \sum_{l=1}^N |m_{\rho(i)\,\rho(l)}| 
                = ||M^\rho||_\infty \nonumber
  \end{eqnarray}
For $||M||_{\infty}$ and $||M||_{\max}$, the proof is similar.
\item The proofs are similar.
  \begin{eqnarray}
  ||v||_1\,    &=& \sum_{j=1}^N |v_j|           = \sum_{j=1}^N |v_{\rho(j)}|          = ||v^\rho||_1 \nonumber\\\nonumber
  ||v||_\infty &=& \max_{j\in C} |v_j|      = \max_{j \in C} |v_{\rho(j)}|    = ||v^\rho||_\infty \\\nonumber
  ||v||_2\,    &=& \sqrt{\sum_{j=1}^N (v_j)^2 } = \sqrt{\sum_{j=1}^N (v_{\rho(j)})^2} = ||v^\rho||_2
  \end{eqnarray}
\end{enumerate}
\end{proof}

In addition, some properties specific for vectors are given:
\begin{corollary}
Let $C$ be a finite state space with $|C|=N$, $\rho$ be a permutation of $C$ (with the \mbox{permutation} matrix $T$), 
$v^1,\,v^2$ $N$-dimensional vectors.\\
Then the following properties hold:
\begin{enumerate}
\item Scalar product: $<v,w>=<v^\rho,w^\rho>$.
\item Outer product: $M=v\otimes w$ iff $M^\rho=v^\rho\otimes w^\rho$.
\end{enumerate}
\end{corollary}
This can be proved based on lemma~\ref{lemreorder}.

\begin{corollary}[Reordering of the exponential matrix]\label{symexpmat} 
For any matrix $M$ it holds that $\mbox{\rm Exp}(M^\rho)=\left(\mbox{\rm Exp}(M)\right)^\rho$. 
\end{corollary}
\begin{proof}
The exponential of a matrix is defined as $\mbox{Exp}(M)=\sum_i^\infty \frac{M^i}{i!}$.\\
The reordering holds for all summands, hence also for the sum.
\end{proof}

\subsection{Basis of vector spaces} \label{basisvecex}

This section shall present the basis of the vector spaces of state symmetric matrices and vectors for all 
non-isomorphic subgroups of state symmetries of sets with maximally four elements. 
To describe the basis we use an abreviated notation of vectors and matrices which define
groups of pure state symmetries (compare definition~\ref{puress} and theorem~\ref{charvecgroup}).
The basis can be obtained by taking the vectors and matrices to pieces by taking 
one elements as 1 and the other to 0.

\subsubsection{Sets of maximally three elements}

$N=|C|$

~

\begin{center}
\begin{tabular}{||c||c|r||c|c|c||c|c||} \hline\hline
   & $N$ & $N^2$ & $\GG$ & abstract group & |\GG| & $v$ & $M$ \\ \hline\hline
a) & 1 & 1 & $\{id\}$ & $\CyG_1$ & 1
       & $(1)$          
       & $\begin{pmatrix} 1 \end{pmatrix}$ \\ \hline
b) & 2 & 4 & $<(1\,2)>$ & $\CyG_2$ & 2
       & $\begin{pmatrix} 1 \\ 1 \end{pmatrix}$
       & $\begin{pmatrix} 1 & 2 \\ 2 & 1 \end{pmatrix}$ \\ \hline
c) & 3 & 9 & $<(1\,2)(3)>$ & $\CyG_2$ & 2
       & $\begin{pmatrix} 1 \\ 1 \\ 2  \end{pmatrix}$
       & $\begin{pmatrix} 1 & 2 & 3 \\ 2 & 1 & 3 \\ 4 & 4 & 5 \end{pmatrix}$ \\ 
   & & & & & & & \\
d) & & & $<(1\,2\,3)>$ & $\CyG_3=\AlG_3$ & 3
       & $\begin{pmatrix} 1 \\ 1 \\ 1 \end{pmatrix}$    
       & $\begin{pmatrix} 1 & 2 & 3 \\ 3 & 1 & 2 \\ 2 & 3 & 1 \end{pmatrix}$ \\ 
   & & & & & & & \\
e) & & & $\begin{array}{c}<(1\,2)(3),\mbox{~ ~}\, \\ \mbox{\,~ ~}(1\,2\,3)>\end{array}$ & $\SyG_3$ & 6
       & $\begin{pmatrix} 1 \\ 1 \\ 1 \end{pmatrix}$    
       & $\begin{pmatrix} 1 & 2 & 2 \\ 2 & 1 & 2 \\ 2 & 2 & 1 \end{pmatrix}$ \\ 
\hline\hline
\end{tabular}
\end{center}

It can be noted, that $\SSy(\maV_{{\cal Z}_3})=\SyG_3$, while $\SSy(\maM_{{\cal Z}_3})=\CyG_3$.
Hence, $\CyG_3$ is a group of pure matrix state symmetries, but not of pure vector state symmetries.

~

In addition, $\maM_{{\cal S}_3}$ consists of only symmetric matrices $M=M\Transp$.
For $\maM_{{\cal Z}_2}$, there is a subspace of symmetric matrices with basis 
$\begin{pmatrix} 1 & 2 & 3 \\ 2 & 1 & 3 \\ 3 & 3 & 5 \end{pmatrix}$, 
while the subspace of symmetric matrices of $\maM_{{\cal Z}_3}$ is equal to $\maM_{{\cal S}_3}$.

\subsubsection{Sets of four elements}

$N=4$, $N^2=16$

~

\begin{center}
\begin{tabular}{||c||c|c|c||c|c||} \hline\hline
  & $\GG$ & abstract group & |\GG| & $v$ & $M$ \\ \hline\hline
a) & ~~~$<(1\,2)(3)(4)>$~~~ & $\CyG_2$ & 2 
       & $\begin{pmatrix} 1 \\ 1 \\ 2 \\ 3  \end{pmatrix}$
       & $\begin{pmatrix} 1 & 2 & 3 & 4 \\ 2 & 1 & 3 & 4 \\ 5 & 5 & 6 & 7 \\ 8 & 8 & 9 & 10 \end{pmatrix}$ \\ 
 & & & & & \\
b) & $<(1\,2)(3\,4)>$ & $\CyG_2$ & 2 
       & $\begin{pmatrix} 1 \\ 1 \\ 2 \\ 2  \end{pmatrix}$
       & $\begin{pmatrix} 1 & 2 & 3 & 4 \\ 2 & 1 & 4 & 3 \\ 5 & 6 & 7 & 8 \\ 6 & 5 & 8 & 7  \end{pmatrix}$ \\ 
 & & & & & \\
c) & $\begin{array}{c}<(1\,2)(3)(4),\mbox{~ ~}\, \\ \mbox{\,~ ~}(1)(2)(3\,4)>\end{array}$ & $\CyG_2 \times \CyG_2$ & 4
       & $\begin{pmatrix} 1 \\ 1 \\ 2 \\ 2  \end{pmatrix}$
       & $\begin{pmatrix} 1 & 2 & 3 & 3 \\ 2 & 1 & 3 & 3 \\ 4 & 4 & 5 & 6 \\ 4 & 4 & 6 & 5  \end{pmatrix}$ \\ 
 & & & & & \\
d) & $\begin{array}{c}<(1\,2)(3\,4), \mbox{~ ~} \\ \mbox{\,~ ~}(1\,3)(2\,4)> \end{array} $ & $\CyG_2 \times \CyG_2$ & 4
       & $\begin{pmatrix} 1 \\ 1 \\ 1 \\ 1  \end{pmatrix}$
       & $\begin{pmatrix} 1 & 2 & 3 & 4 \\ 2 & 1 & 4 & 3 \\ 3 & 4 & 1 & 2 \\ 4 & 3 & 2 & 1  \end{pmatrix}$ \\ 
 & & & & & \\
e) & $<(1\,2\,3\,4)>$ & $\CyG_4$ & 4       & $\begin{pmatrix} 1 \\ 1 \\ 1 \\ 1  \end{pmatrix}$
       & $\begin{pmatrix} 1 & 2 & 3 & 4 \\ 4 & 1 & 2 & 3 \\ 3 & 4 & 1 & 2 \\ 2 & 3 & 4 & 1  \end{pmatrix}$ \\ 
 & & & & & \\
f) & $\begin{array}{c}<(1\,2\,3\,4), \mbox{~ ~} \\ \mbox{\,~ ~ ~} (1\,2)(3\,4)> \end{array} $ & $\CyG_4\times\CyG_2$ & 8
       & $\begin{pmatrix} 1 \\ 1 \\ 1 \\ 1  \end{pmatrix}$
       & $\begin{pmatrix} 1 & 2 & 3 & 3 \\ 2 & 1 & 3 & 3 \\ 3 & 3 & 1 & 2 \\ 3& 3 & 2 & 1  \end{pmatrix}$ \\ 
\hline\hline
\end{tabular}

\begin{tabular}{||c||c|c|c||c|c||} \hline\hline
 & $\GG$ & abstract group & |\GG| & $v$ & $M$ \\ \hline\hline
g) & $<(1\,2\,3)(4)>$ & $\CyG_3$ & 3
       & $\begin{pmatrix} 1 \\ 1 \\ 1 \\ 2  \end{pmatrix}$
       & ~$\begin{pmatrix} 1 & 2 & 3 & 4 \\ 3 & 1 & 2 & 4 \\ 2 & 3 & 1 & 4 \\ 5 & 5 & 5 & 6  \end{pmatrix}$~ \\ 
 & & & & & \\
h) &  $\begin{array}{c}<(1\,2\,3)(4), \mbox{~ ~} \\ \mbox{\,~ ~ ~} (1\,2)(3)(4)> \end{array}$ & $\SyG_3$ & 6
       & $\begin{pmatrix} 1 \\ 1 \\ 1 \\ 2 \end{pmatrix}$    
       & $\begin{pmatrix} 1 & 2 & 2 & 3 \\ 2 & 1 & 2 & 3 \\ 2 & 2 & 1 & 3 \\ 4 & 4 & 4 & 5 \end{pmatrix}$ \\
 & & & & & \\
i) & $\begin{array}{c}<(1\,2\,3)(4), \mbox{~ ~} \\ \mbox{\,~ ~ ~} (1\,2)(3\,4)> \end{array}$ & $\AlG_4$ & 12
       & $\begin{pmatrix} 1 \\ 1 \\ 1 \\ 1 \end{pmatrix}$    
       & $\begin{pmatrix} 1 & 2 & 2 & 2 \\ 2 & 1 & 2 & 2 \\ 2 & 2 & 1 & 2 \\ 2 & 2 & 2 & 1 \end{pmatrix}$ \\
 & & & & & \\
j) & $\begin{array}{c}<(1\,2\,3\,4), \mbox{~ ~ ~} \\ \mbox{~ ~ ~ ~} (1\,2)(3)(4)> \end{array}$ & $\SyG_4$ & 24
       & $\begin{pmatrix} 1 \\ 1 \\ 1 \\ 1 \end{pmatrix}$    
       & $\begin{pmatrix} 1 & 2 & 2 & 2 \\ 2 & 1 & 2 & 2 \\ 2 & 2 & 1 & 2 \\ 2 & 2 & 2 & 1 \end{pmatrix}$ \\
\hline\hline
\end{tabular}
\end{center}

a) and b) lead to different spaces for matrices and vectors, even if their groups of state symmetries are isomorphic.
The same holds for c) and d).

d) and e) result in different spaces for matrices, because d) has only symmetric matrices and e) has not
(diffent space means: the basis matrices cannot be made equal by renumbering and reordering).

Also f), i), and j) consists of only symmetric matrices.

~

Only a), c), h) and j) describe groups of pure vector symmetries.

Note that $\maM(\SyG_4)=\maM(\AlG_4)$, so that $\AlG_4$ is not a group of pure matrix state symmetries.

\end{appendix}

\section*{Symbols}

\begin{small}
\begin{tabular}{ll}
$C$                                      & finite set (state space $\{1,\ldots,N\}$) \\
$\calC=\{C_i\}_i$                        & partition of $C$ \\
$\Rspace$                                & vector space of reals  \\
$\Mspace$                                & vector space of real valued $N\times N$ matrices \\
$GL_N({\mathbb R})$                      & group of invertible square matrices \\[2mm]

$v$;         $L_v$                       & column vector; list of column vectors \\
$M$;         $L_M$                       & matrix; list of matrices \\
$\rho$;      $P$;        $\GG$           & permutation (state symmetry); list of permutations; group of permutations \\
$T$;         $T(\rho)$                   & permutation matrix; permutation matrix of $\rho$ \\
$<\rho>$                                 & group, generated by $\rho$\\
$v^\rho$;    $M^\rho$                    & reordered column vector; reordered matrix \\
$M\Transp$                               & Transposition of $M$ \\[2mm]

$\maV$;      $\maV(L_v)$                 & vector space of column vectors; vector space, spanned by elements of $L_v$ \\
$\maV_\rho$; $\maV_P$;   $\maV\sGG$      & vector spaces of state symmetric column vectors \\
$\maM$;      $\maM(L_M)$                 & vector space of matrices; vector space, spanned by elements of $L_M$ \\
$\maM_\rho$; $\maM_P$;   $\maM\sGG$      & vector spaces of state symmetric matrices\\[2mm]

$\SSy(.)$                                & group of state symmetries of object .\\
$\SyG_N$;    $\AlG_N$;   $\CyG_N$        & Symmetric group; Alternate group; Cyclic group of order $N$\\
$\SyG(C)$                                & Symmetric group of the elements of $C$\\[2mm]

$u=(v^{(k)})_{k=0,\,1\,\ldots}$          & sequence of column vectors \\
$u(i)$                                   & first $i+1$ elements of $u$\\
$\mathbb{W}_{u}$                         & set of generating matrices of $u$\\
$\mathbb{W}_{u}^{\rho}$                  & set of state symmetric generating matrices of $u$\\

\end{tabular}
\end{small}

\section*{Acknowledgements}

The author would like to thank Dr. K.~Richter (University of Halle-Wittenberg)
and Prof. C.~L\"ofvall (Stockholms University) 
for their support during the writing of the underlying 
Master's thesis \cite{Ring:96} and continuous helpful discussions on this topic.


\end{document}